\documentclass{article}

\usepackage{arxiv}

\usepackage[utf8]{inputenc} 
\usepackage[T1]{fontenc}    
\usepackage{hyperref}       
\usepackage{url}            
\usepackage{booktabs}       
\usepackage{amsfonts}       
\usepackage{nicefrac}       
\usepackage{microtype}      
\usepackage{lipsum}

\usepackage{amsmath,amsthm,amssymb,graphicx,arxiv,longtable,booktabs,array,siunitx,float,placeins}

\usepackage[aboveskip=0pt,belowskip=0pt]{subcaption}
\usepackage[tableposition=below]{caption}
\usepackage{multicol}
\usepackage{color}
\usepackage{url}
\providecommand{\Sk}{\mathrm{Sk}}
\providecommand{\Ku}{\mathrm{Ku}}

\newcommand{\bo}{\mathbf}
\newcommand{\bs}{\boldsymbol}
\newtheorem{The}{Theorem}[section]
\newtheorem{Def}{Definition}[section]

\numberwithin{equation}{section}

\title{Notion of information and independent component analysis}
\date{January 29, 2020}
\author{Una Radojicic  \\
  Vienna University of Technology\\  
  \texttt{una.radojicic@tuwien.ac.at} \\
  \And
        Klaus Nordhausen \\
  Vienna University of Technology\\
  \texttt{klaus.nordhausen@tuwien.ac.at} \\
  \And
        Hannu Oja \\
  University of Turku\\
  \texttt{hannu.oja@utu.fit} \\
}

\begin{document}
\maketitle

\begin{abstract}
   Partial orderings and measures of information for continuous univariate random variables
   with special roles of Gaussian and uniform distributions are discussed.
   The information measures and measures of non-Gaussianity including third and fourth cumulants are generally used as projection indices in the projection pursuit approach for the independent component analysis. The connections between information,  non-Gaussianity and  statistical independence in the context of independent component analysis is discussed in detail.
\end{abstract}

\keywords{Dispersion\and entropy\and kurtosis\and partial orderings}

\section{Introduction}\label{sec1}

In the engineering literature independent component analysis (ICA) \cite{hy,IC2} is often described as a search for
the uncorrelated linear combinations of the original variables that maximize non-Gaussianity.
The estimation procedure then usually has two steps.  First, the vector of principal components is found and the components are standardized to have zero means and unit variances, and second, the vector is further rotated so that the new components  maximize a selected measure of non-Gaussianity. It is then argued that the components obtained in this way are made as independent as possible or that they display the components with maximal information. \cite{hy} for example give a heuristic argument that, according to the central limit theorem, weighted sums of independent non-Gaussian random variables are closer to Gaussian than the original ones. In this paper, we discuss and clarify the somewhat vague connections between non-Gaussianity, independence and  notions of information in the context of the independent component analysis.

In Section~\ref{sec2} we first introduce descriptive measures for location, dispersion, skewness and kurtosis of univariate random variables with some discussion of corresponding partial orderings.
In this part of the paper we assume that the considered  univariate random variable $x$ has a finite mean $E(x)$ and  variance $Var(x)$, cumulative distribution function  $F$  and continuously differentiable probability density function  $f$. Skewness, kurtosis and other  cumulants of the standardized variable $(x-E(x))/\sqrt{Var(x)}$ are often used to measure non-Gaussianity of the distribution of $x$.  The most popular measures of statistical information are the differential entropy $H(f)=-\int f(x) \log(f(x))dx$ and
the Fisher information in the location model, that is, $J(f)=\int f(x) [f'(x)/f(x)]^2 dx$.
These and other information measures with related partial orderings and their use as measures of non-Gaussianity are discussed in the later part of Section~\ref{sec2}.

The multivariate independent components model is discussed  in Section~\ref{sec3}.  It is then
assumed that, for a $p$-variate random vector $\bo x$, there is a linear operator $\bo A\in \mathbb{R}^{p\times p}$ such that
$\bo A\bo x$ has independent components. Under certain assumptions, the projection pursuit approach can be used to find the rows of $\bo A$ one-by-one and various information measures as well as cumulants have been used as projection indices. In Section~\ref{sec3} the connections between non-Gaussianity, independence and information in this context is discussed in detail. The paper ends with some final remarks in Section~\ref{sec4}.

\section{Some characteristics of a univariate distribution}\label{sec2}

\subsection{Location, dispersion, skewness and kurtosis}\label{subsec2.1}

We consider a continuous
random variable $x$ with the finite mean $E(x)$, finite variance $Var(x)$, density function $f$ and
cumulative density function  $F$.
Location, dispersion, skewness and kurtosis are often considered by defining
the corresponding measures or functionals for these properties.
Location and dispersion measures, write $T(x)$ and $S(x)$,
are functions of the distribution of $x$ and defined as follows.

\begin{Def}\label{location_and_dispersion}
\ \
\begin{enumerate}
\item $T(x)\in \mathbb{R}$ is a location measure if $T({ax+b})=aT(x)+b$, for all
$\forall a,b\in \mathbb{R}$.
\item $S(x)\in \mathbb{R}_{+}$ is a dispersion measure if $S({ax+b})=|a|S(x)$, for all
$\forall a,b\in \mathbb{R}$.
\end{enumerate}
\end{Def}

Clearly, if $T$ is a location measure and $x$ is symmetric around $\mu$, then $T(x)=\mu$ for all location measures. For squared dispersion measures $S^2$,
\cite{huber} considered the concepts of additivity,  subadditivity and superadditivity.
These concepts appear to be crucial in developing tools for the independent component
analysis and are defined as follows.

\begin{Def}\label{additivity}
Let $S^2$ be  a squared dispersion measure.
\begin{enumerate}
\item $S^2$ is {\it additive}  if
$S^2(x+y)=  S^2(x)+S^2(y)$ for all independent $x$ and $y$.
\item $S^2$ is {\it subadditive}  if
$S^2(x+y)\le  S^2(x)+S^2(y)$ for all independent $x$ and $y$.
\item $S^2$ is {\it superadditive} if
 $S^2(x+y)\ge S^2(x)+S^2(y)$ for all independent $x$ and $y$.
\end{enumerate}
\end{Def}

The mean $E(x)$ and the variance $Var(x)$ are important and most popular  location and squared dispersion measures.
It is well known that   $Var(x+y)=Var(x)+Var(y )$ for independent $x$ and $y$, and $E(x+y)=E(x)+E(y)$ is true even for dependent $x$ and $y$.
These additivity properties are highly important in certain applications and in fact characterize the mean and variance among continuous measures as follows.

\begin{The}\label{additivity theorem}
\ \
\begin{enumerate}
\item Let a location measure $T$ be additive and  continuous at $N(0,1)$, that is,  $z_n\to_d z\sim N(0,1)$ implies that $T(z_n)\to T(z)=0$. Then $T(x)=E(x)$ for all $x$ with finite second moments.
\item Let a squared dispersion measure $S^2$ be additive and continuous at $N(0,1)$, that is, $z_n\to_d z\sim N(0,1)$ implies that $S^2(z_n)\to S^2(z)>0$. Then $S^2(x)=S^2(z) \ Var(x)$ for all $x$ with finite second moments.
\end{enumerate}
\end{The}

Comparison of different location measures $T_1$ and $T_2$ and dispersion measures $S_1$ and $S_2$,
provides measures of skewness and kurtosis as
$$  \Sk(x)=\frac{T_1(x)-T_2(x)}{S(x)} \ \ \mbox{and}\ \
    \Ku(x)= \frac {S_1^2(x)}{S_2^2(x)}.
$$
Classical measures of skewness and kurtosis proposed in the literature can be written in this way.
Note that
both measures are affine invariant in the sense that
$$
\Sk({ax+b})= sgn(a) \Sk(x) \ \ \mbox{and}\ \
\Ku({ax+b})=  \Ku(x).
$$
If $x$ has a symmetric distribution, then $\Sk(x)=0$. In the literature, kurtosis measures are thought to measure
the peakedness and/or the heaviness of the tails of the density of $x$ but, as we will see
in Section~\ref{subsec2.3}, $\Ku(x)$ as defined here
may be a global measure of deviation from the normality
and have also been used  as an affine invariant information measure for some special choices of the dispersion measures $S_1$ and $S_2$.

Moment and cumulant generating functions defined as
\[
E\left[e^{tx} \right]=\sum_{k=0}^\infty {\mu_kt^k}{/k!} \ \
\mbox{and}\ \
\log E\left[e^{tx} \right]=\sum_{k=0}^\infty  {\kappa_kt^k}{/k!}
\]
respectively, generate classical measures, i.e.,
moments
$E(x)=\mu_1(x)$ and  $Var(x)=\mu_2(x-\mu_1(x))$ and cumulants $\kappa_3(x^{st})$ and
$\kappa_4(x^{st})$ where  $x^{st}={(x-E(x))}/{\sqrt{Var(x)}}$.
The cumulants $\kappa_k$, $k=1,2,...$  are additive as $\log E\left[e^{tx} \right]$ is additive,
and $\kappa_k^{2/k}(x-E(x))$, $k=2,3,...$ are subadditive squared dispersion measures which follows from the Minkowski inequality, see  \cite{huber}.
Another class of measures is given by the quantiles
$q_u=F^{-1}(u)$, $0<u<1$, with corresponding measures such as
\[
q_{1/2},\ \   q_{1-u}-q_{u},\ \  \frac {q_u+q_{1-u}-2q_{1/2}}{q_{1-u}-q_{u}}, \  \ \mbox{and} \ \ \frac{q_{1-u}-q_{u}} {q_{1-v}-q_{v}},\ \ \
0<u<v<\frac 12.
\]
These quantile based measures provide robust alternatives to moment based measures.
To our knowledge, they however lack
the additivity properties stated in Definition~\ref{additivity} which makes them unsuitable for usage in the independent component analysis.

An alternative strategy to consider the properties of distributions is to define {\it partial orderings}
for location, dispersion, skewness and kurtosis.
For continuous $x$ and $y$ with cumulative distribution functions  $F$ and $G$,
write $\Delta(x)=G^{-1}(F(x))-x$. The function $\Delta(x)$ is called a shift function of
$x$ as $x$ when shifted by $\Delta(x)$  and has the distribution of $y$.
The transformation $x\mapsto x+\Delta(x)$ is also known as the (univariate) Monge-Kantorovich optimal transport map.
Using function $\Delta$ we can naturally define the following partial orderings \cite{Lehman1,Lehman2,Zwet,Oja91}.
\begin{multicols}{2}
\begin{enumerate}
\item Location ordering: $\Delta$ is positive.
\item Dispersion ordering: $\Delta$ is increasing.
\item Skewness ordering: $\Delta$ is convex.
\item Kurtosis ordering: $\Delta$ is concave-convex.
\end{enumerate}
\end{multicols}
\cite{Lehman1,Lehman2,Oja91} then stated that, in addition to the affine equivariance and invariance properties,  the measures of location, dispersion, skewness and kurtosis should be monotone with respect to corresponding orderings. For finding monotone measures in the dispersion case, for example, $\Delta$ is increasing if and only if
\[
\mbox{$E[C(x-E(x))]\le E[C(y-E(y))]$ for all convex $C$.}\
\]
which is also called the dilation order. It implies for example that the  measures $(E[ |x-E(x)|^k] )^{1/k}$, $k>1$,  are monotone dispersion measures.

\subsection{Information and discrete distributions}\label{subsec2.2}

Consider a discrete random variable with $k$ possible values (`alphabets')  with probabilities listed in
$p=(p_1,...,p_k)$. Write $p_{(1)}\le ... \le p_{(k)}$ for the ordered probabilities.  It is sometimes presumed that a distribution $p$ is informative if it can provide `surprises' with very small $p_i$'s.  On the other hand, people often  claim that $p$ is informative if
the result of the experiment is known with a high probability, that is, if  only  one or few values have high $p_i$'s. These somewhat naive characterizations  suggest the following well-known partial ordering  for
discrete distributions \cite{Marshall}.

\begin{Def}
{\it Majorization} :
$p \prec q$ if $\displaystyle \sum_{i=1}^j p_{(i)}\ge \sum_{i=1}^j q_{(i)}$, $j=1,...,k$, and then $p$ is said to be majorized by $q$.
\end{Def}

Majorization is nothing but a dispersion ordering (and a dilation order) for the discrete distributions
with $k$ equiprobable values $p_1,...,p_k$ in $[0,1]$ with mean $1/k$.
Then, according to \cite{po},
\begin{eqnarray*}
\mbox{$p \prec q$} \ &\Leftrightarrow &
\mbox{$p=q \bo L$ with some doubly stochastic matrix $\bo L$}  \ \ \\ &\Leftrightarrow &
\mbox{$\sum_{i=1}^k C(p_i) \le \sum_{i=1}^k C(q_i) $ for all continuous convex $C$.}
\end{eqnarray*}
\
The doubly stochastic matrix $\bo L$ is a matrix with non-negative elements such that all row sums and all column sums are one. The doubly stochastic operator $\bo L$ is then in fact a convex combination of permutations; \ $p$ is obtained from $q$ by this `smoothing' and is therefore less informative.
Further, for all $p$,  $$(1/k,...,1/k)\prec p \prec (0,...,0,1)$$  and,
for simple mixtures,
$p \prec q \ \ \Rightarrow \ \  p \prec \lambda p + (1-\lambda) q \prec q, \ 0\le \lambda \le 1.$

We can now give the following.
\begin{Def}
Let $p=(p_1,...,p_k)$ list the probabilities of $k$ possible values of a discrete random variable, that is, $p_1,...,p_k\in [0,1], \sum_{i=1}^k p_i=1$.
A measure $M(p)$ is a information measure if it is monotone with respect to majorization.
\end{Def}

Note that, as $(p_1,...,p_k)\prec (p_{(1)},...,p_{(k)}) \prec (p_1,...,p_k)$, the definition implies that the information measures are invariant under permutations
of the probabilities in $(p_1,...,p_k)$.
The equivalent conditions for majorization then suggest quantities such as
\[
H(p) = -\sum_{i=1}^k \log (p_i) p_i, \ \ \ H^*(p)=\sum_{i=1}^k p_i^2 \ \ \ \mbox{and}\ \
\ \ H^{**}(p)=p_{(k)}
\]
and $-H$, $H^*$ and $H^{**}$ are monotone information measures that easily extend to continuous and multivariate cases. The {\it Shannon's entropy} \cite{Shannon} $-\sum_{i=1}^k \log_2 (p_i) p_i$ is often
seen as a measure of  ability to compress the data \ (e.g. lower bound for the  expected number of bits to store the data).

\subsection{Some information measures for continuous distributions}\label{subsec2.3}

Consider next a continuous random variable $x$ with the continuously differentiable probability density function $f$ and finite variance $Var(x)$.
The three measures from the discrete case straightforwardly extend in the continuous case to
\begin{eqnarray*}
  H(x) &=& -E[\log f(x)]=-\int_{-\infty}^{\infty} f(x) \log f(x) dx, \\
  H^*(x) &=& E[f(x)]= \int_{-\infty}^{\infty} f^2(x) dx, \ \ \mbox{and} \\
  H^{**}(x) &=& \sup_x f(x) = f(x_{mode}),\ \ \mbox{if the mode $x_{mode}$ exists}.
\end{eqnarray*}
The Fisher information in the location model $f(\cdot-\mu)$ at $\mu=0$ given by
\[
J(x)=\int_{-\infty}^{\infty} f(x)\left( \frac {f'(x)}{f(x)} \right)^2 dx.
\]
is also often used as an information measure \cite{kb}.

 The measure $H(x)$ is popular in the literature and known as the {\it differential entropy}. Under certain restrictions, the measure has the following maximizers \cite{it}. For the distributions on $\mathbb{R}$ with a fixed variance, $H(x)$ is maximized
if $x$ has a normal distribution. For distributions on $ \mathbb{R}_{+}$ with a fixed mean, $H(x)$
is maximized at the exponential distribution. For distributions on a finite interval, $H(x)$ is maximized at the uniform distribution on that interval. Note  that, in the Bayesian analysis,  these three distributions
are often used as priors that reflect `total ignorance'.

We next show that the three straightforward extensions $H$, $H^*$ and $H^{**}$ as well as
the Fisher information $J$ provide squared dispersion measures as in Definition~\ref{location_and_dispersion} but with an interesting additional invariance property. First note that the measures are invariant under location shift of the distribution but not
under rescaling of the variable.
Recall that information as stated for discrete distributions is invariant under the permutations of
the probabilities in $(p_1,...,p_k)$. All permutations consist of successive pairwise exchanges
 of two probabilities.  In the continuous case,  similar elemental probability density transformations may be  constructed as follows.  For all  $a<a+\Delta<b<b+\Delta$ and density function $f$, write
$$
 f_{a,b,\Delta}(x) =\left\{
  \begin{array}{ll}
    f(x), & \hbox{$x\in \mathbb{R}-[a,a+\Delta]-[b,b+\Delta] $} \\
    f(b+(x-a)), & \hbox{$x\in[a,a+\Delta] $} \\
    f(a+(x-b)), & \hbox{$x\in[b,b+\Delta] $}
  \end{array}
\right.
$$
The transformation allows the manipulation of the properties of the distribution in many ways. The transformation can for example be used to move some probability mass from the centre of distribution to the tails and in this way to  manipulate the variance and the kurtosis of the distribution
for example. As far as we know, this transformation has not been discussed in the literature.
It is surprising that the information  measures $H$, $H^*$, $H^{**}$ and $J$  provide  dispersion measures which are  invariant under these transformations.

\begin{The}\label{dispersion measures}
The {\it entropy power} $e^{2H(x)}$ and measures  $[H^*(x)]^{-2}$, $[H^{**}(x)]^{-2}$
and $[J(x)]^{-1}$  are
squared dispersion measures that are  invariant under the transformations $f\to  f_{a,b,\Delta}$.
The measures $e^{2H(x)}$ and $[J(x)]^{-1}$ are superadditive.
\end{The}

\subsection{Affine invariant information measures}\label{subsec2.3b}

We now further discuss the properties of the dispersion measures in Theorem~\ref{dispersion measures}
and, to find affine invariant information measures, consider  the ratios of the variance to these squared dispersion measures. The  ratio of the variance to the entropy power, that is,  $Var(x)e^{-2H(x)}$  is minimized at the normal distribution \cite{it}. In a neighbourhood of a normal distribution the negative entropy $-H(x)$ has an interesting approximation using third and fourth cumulants. \cite{Sibons} showed that the negative differential entropy for the density $f(x)=\varphi(x)(1+\epsilon(x))$ where $\varphi$ is the density of $N(0,1)$ and $\epsilon$ is a well-behaved ``small'' function that satisfies $E[\epsilon(z)z^k]=0$, $z\sim N(0,1)$,  $k=0,1,2$, can be approximated by $(1/2)\int \varphi(x)\epsilon^2(x)dx \approx
(\kappa_3^2(x)+(1/4) \kappa_4^2(x))/12$. \

Next, $[H^*(x)]^{-2}$  is a (squared) dispersion measure,  and  therefore \ $[H^*(x)]^{2} Var(x)$ provides an affine invariant information  measure. For symmetric distributions, it preserves the concave-convex kurtosis ordering of van Zwet and
$12 [H^*(x)]^{2}  Var(x)$
is in fact the efficiency of the Wilcoxon rank test with respect to  the $t$-test.
Also, for symmetric distributions,  $4 [H^{**}(x)]^{2}  Var(x)$  is a kurtosis measure in  the van Zwet sense and simultaneously the efficiency of the sign test with respect to  the $t$-test.  We also
mention that, if
$ Q(x)= E\left[f(F^{-1}(u))/\varphi(\Phi^{-1}(u))\right]$
with $u\sim U(0,1)$,
then $Q^{-2}(x)$ is a squared dispersion measure and $Q^2(x)Var(x)$ is the efficiency of
the van der Waerden test with respect to the $t$-test in the symmetric case. By the Chernoff-Savage theorem, it attains its minimum 1 at the normal distribution. See \cite{Savage,Lehman3}.

Finally,
the information measure  $Var(x)J(x)\ge 1$ is  minimized
at the normal distribution.  In the location estimation problem in the symmetric case, $Var(x)J(x)$
is also the asymptotic relative efficiency  of the maximum likelihood estimate of the symmetry centre with respect to the sample mean \cite{Serfling}.

\subsection{Information orders for continuous distributions}\label{subsec2.4}

\begin{figure}[h!]
\includegraphics[width=0.95\textwidth,height=0.95\textheight,keepaspectratio]{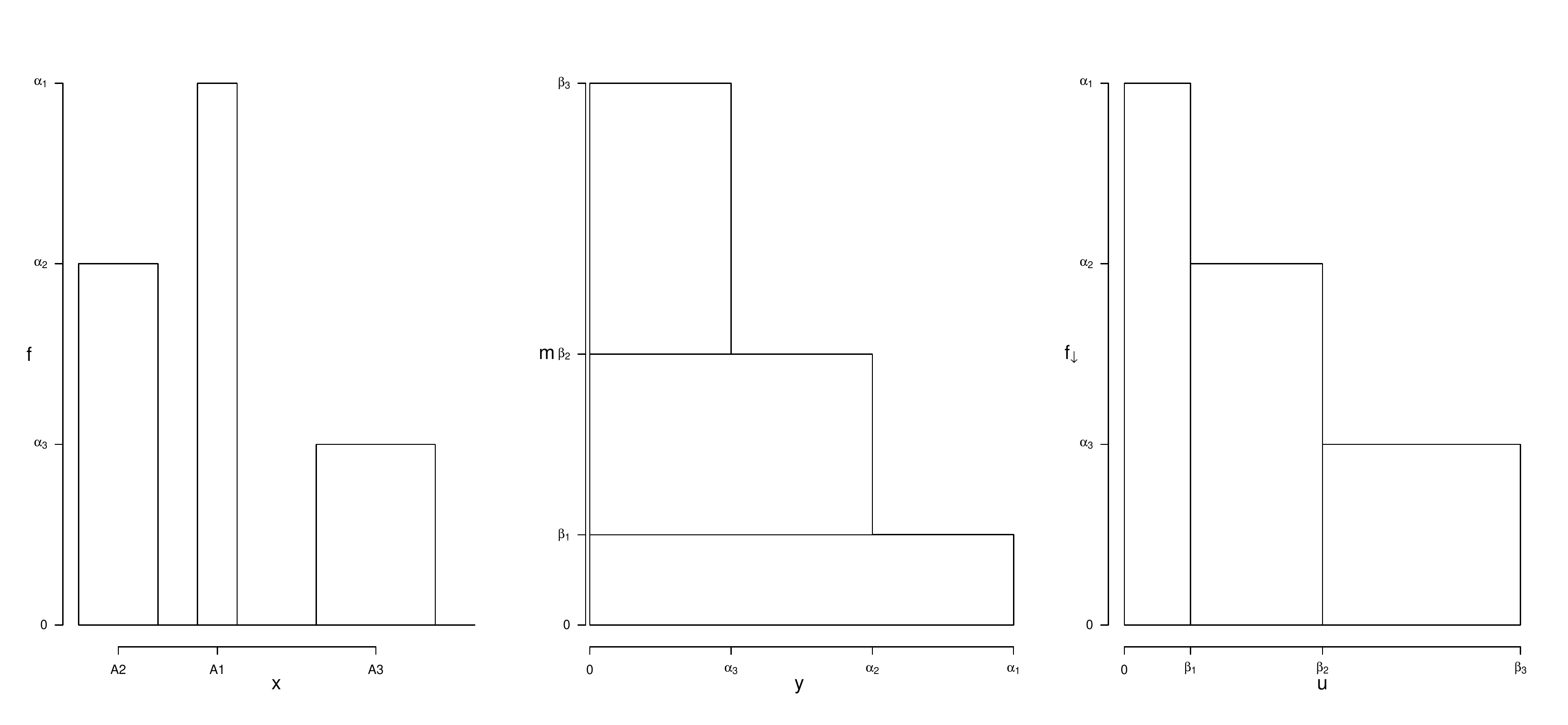}
\caption{Simple function $f$ (left), its distribution function $m$ (middle) and decreasing rearrangement $f_\downarrow$ (right)}
\label{fig:simple function}
\end{figure}

We next outline how to construct partial orderings for information in the univariate continuous case. Let first $x$ be a continuous random variable  with density $f$ on $(0,1)$.
If $m(y)=\mu\{u:f(u)>y\}$ where $\mu$ is Lebesgue measure, then the function
$
f_\downarrow(u)=\sup\{y:m(y)>u\},\ \ u\in (0,1),
$
provides  the {\it decreasing rearrangement} of $f$. Note that any density function on $(0,1)$ can be approximated by a simple density function
$f(x)=\sum_{i=1}^{k}\alpha_i\chi_{A_j}(x)$, where $\alpha_1<\alpha_2<\cdots<\alpha_k$, and
$A_1,...,A_k$ are disjoint Lebesque-measurable sets on $(0,1)$ and $\chi_A$ is the characteristic function of set $A$. Then
$$
       m(y)=\sum_{i=1}^{k}\beta_i\chi_{B_i}(y) \ \ \mbox{and}\ \ f_\downarrow(u)=\sum_{i=1}^{k}\alpha_i\chi_{[\beta_{i-1},\beta_i)}(u),
       $$
where $\beta_i=\sum_{j=1}^{i}\mu(A_j)$, $B_i=[\alpha_{i+1},\alpha_i)$ for $i=1,2,\ldots,k$, and $\beta_0=\alpha_{k+1}=0$.
For a better insight, see Figure~ \ref{fig:simple function}.
For more details and examples, see e.g. \cite{Kristiansson}.

Using the decreasing rearrangement we can give the following definitions.

\begin{Def}\label{partial}
Let $f$ and $g$ be density functions on the interval $(0,1)$. Then $g$
has more information than $f$, write  $f\prec g$,  if
\[
\int_0^u f_\downarrow(v) dv \le \int_0^u g_\downarrow(v) dv,\ \ \mbox{for all $u\in (0,1)$}
\]
\end{Def}

\begin{Def}\label{info_cont}
Let $\mathcal{F}_{(0,1)}$ be the set of density functions $f$ on the interval $(0,1)$.
Then $M_{(0,1)}:\mathcal{F}_{(0,1)}\to \mathbb{R}$ is an information measure if it is monotone
with respect to the partial ordering in Definition~\ref{partial}
\end{Def}
The distribution with minimum information is the uniform distribution on $(0,1)$.
Information measures are
easily found, see \cite{Ryff}, as  $f\prec g$ if and only if
\[
 \int_0^1 C(f(u))du \le \int_0^1 C(g(u))du
\ \ \mbox{for all continuous convex functions $C$}
\]
\cite{Ryff}
also discusses how to construct
linear operators $L$ for which
 $f=Lg \prec g$ when $f\prec g$. \

Consider next  a continuous random variable $x$ on $\mathbb{R}$ with pdf $f$.
To find a location and scale-free version of the density, \cite{Staudte} 
proposed the transformation
\[
f(x),\ x \in \mathbb{R} \ \ \to \ \ f^*(u)=\frac {f(F^{-1}(u))}{ H^*(x)}, \ u\in (0,1).
\]
Then $f^*$, called {\it the probability density quantile (pdQ)}, is a probability density function on $(0,1)$ which is invariant under linear transformations of the original variable $x$ \cite{Staudte}. It is also true that, for given  $f^*$, the original $f$ is known up location and scale.
Using this density transformation, the definition of an invariant information measure for densities on $\mathbb{R}$ can be given as follows.

\begin{Def}\label{info_cont2}
Let $\mathcal{F}_{\mathbb{R}}$ be a set of density functions $f$ on $\mathbb{R}$
and let $M_{(0,1)}:\mathcal{F}_{(0,1)}\to \mathbb{R}$ be an information measure
for distributions on $(0,1)$. Then $M_{\mathbb{R}}:f\to M_{(0,1)}(f^*)$
is an information measure in the set $\mathcal{F}_{\mathbb{R}}$.
\end{Def}
Note that $M_\mathbb{R}$ is not an extension
of $M_{(0,1)}$ meaning that, $f\in \mathcal{F}_{(0,1)}$ does not imply that  $M_{\mathbb{R}}(f)=M_{(0,1)}(f)$.  $M_\mathbb{R}$ is invariant under rescaling of $f$ while
$M_{(0,1)}$ is not.

Applying Definition~\ref{info_cont2}  and
choosing convex $C(u)=-log(u)$ and $C(u)=log(u) u$,
we get location and scale invariant information measures
for  $f$ such as
\[
\exp\{ -2\int \log(f^*(u)) du \}= e^{2H(x)} [H^*(x)]^2
\]
and
\[
 \exp\{ 4 \int \log(f^*(u)) f^*(u) du \} \ = e^{-2H(f^2/H^*(x))} [H^*(x)]^{-2},
\]
which attain their minimum at the uniform distribution and are invariant under the transformations
$f\to  f_{a,b,\Delta}$. For more details see e.g. \cite{Staudteentropy2018}.

To replace the transformation $f\to f^*$ by a transformation to densities on $(0,1)$ for which minimum information is attained at any density $g$, one can use the following adjustment.

\begin{The} \label{KL}
Let $x$ and $y$ be random variables on $\mathbb{R}$ with the probability density functions $f$ and $g$
and cumulative distribution functions $F$ and $G$, respectively. Then
$$ f:g(u)= \frac {f(G^{-1}(u))} {g(G^{-1}(u))} $$
is a density function on $(0,1)$ and its differential entropy is
$-H(f:g)\ge 0$ is the Kullback-Leibler (KL) divergence  between the distributions
of $x$ and $y$.
\end{The}

Let again $x$ have a density $f$ and let $\varphi$ and $\Phi$ be
the pdf and the cdf of a normal distribution with mean $E(x)$ and variance $Var(x)$.
Then one can show, using similar arguments as in \cite{Staudte}, that
$$ f:\varphi(u)= \frac {f(\Phi^{-1}(u))} {\varphi(\Phi^{-1}(u))},\ \ u\in (0,1) $$
is a location and scale-free density and information measures in  Definition~\ref{info_cont}
applied to the set of densities  $\tilde f=f:\varphi$ attain their
minimums when $f$ has a normal distribution. A collection of information measures is given by
$\int C(\tilde f (u))  du$ with continuous and convex functions $C$ and then we get for example
again
$$
\exp \{2\int log(\tilde f(u))\tilde f (u) du  \}=(2\pi e) e^{-2H(x)}Var(x).
$$

\begin{figure}[h!]
\includegraphics[width=0.95\textwidth,height=0.95\textheight,keepaspectratio]{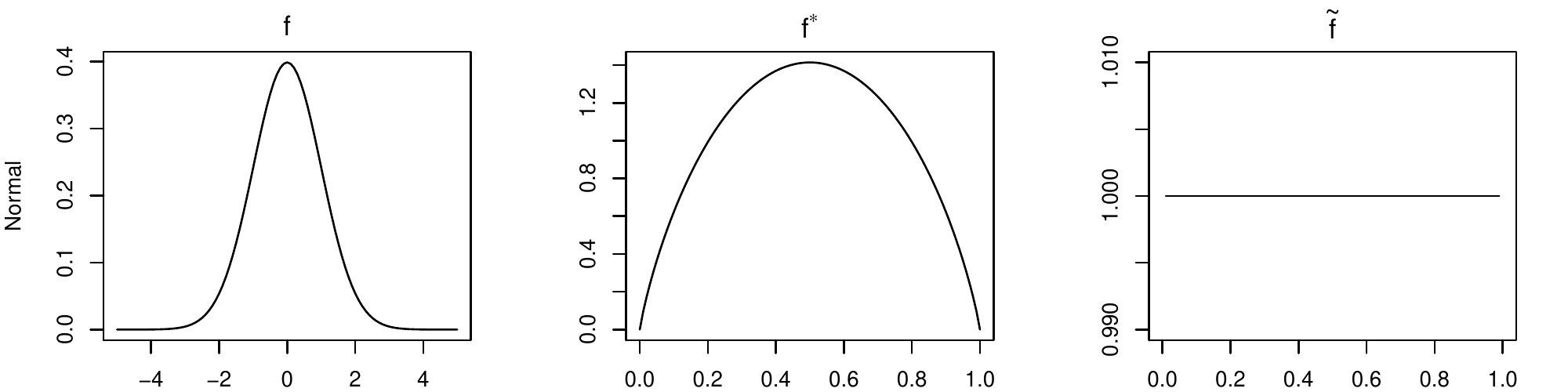}
\includegraphics[width=0.95\textwidth,height=0.95\textheight,keepaspectratio]{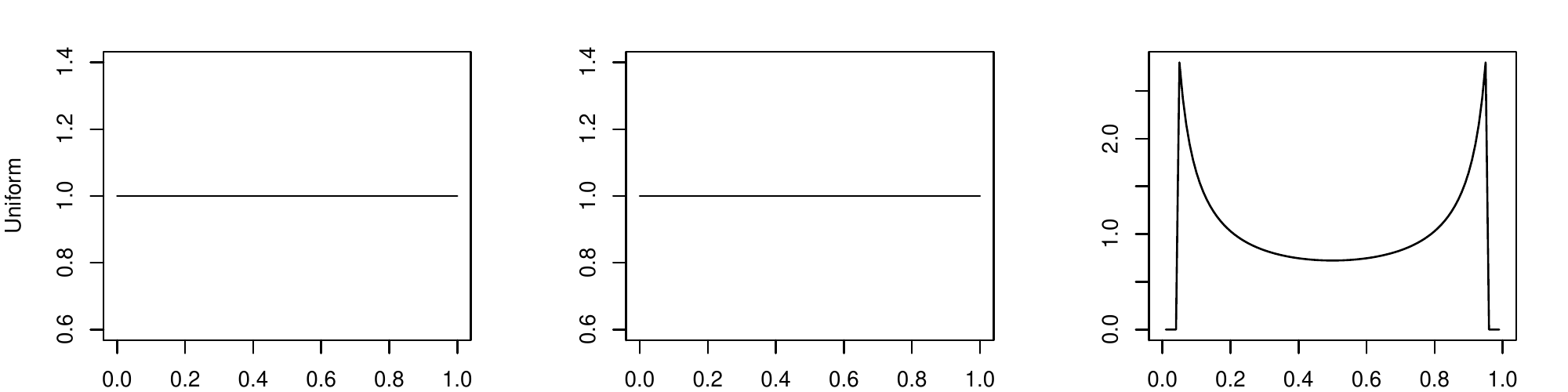}
\includegraphics[width=0.95\textwidth,height=0.95\textheight,keepaspectratio]{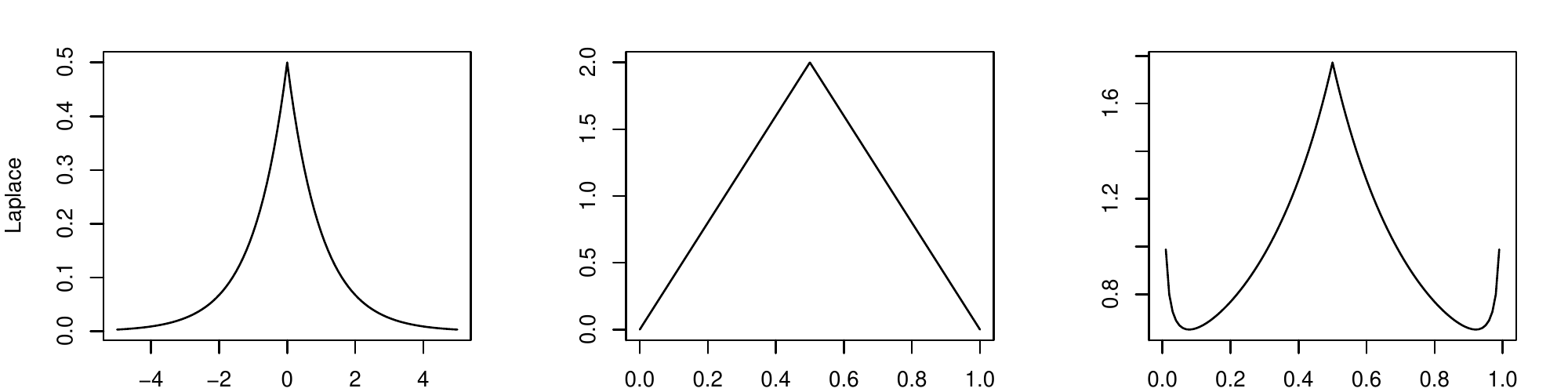}
\includegraphics[width=0.95\textwidth,height=0.95\textheight,keepaspectratio]{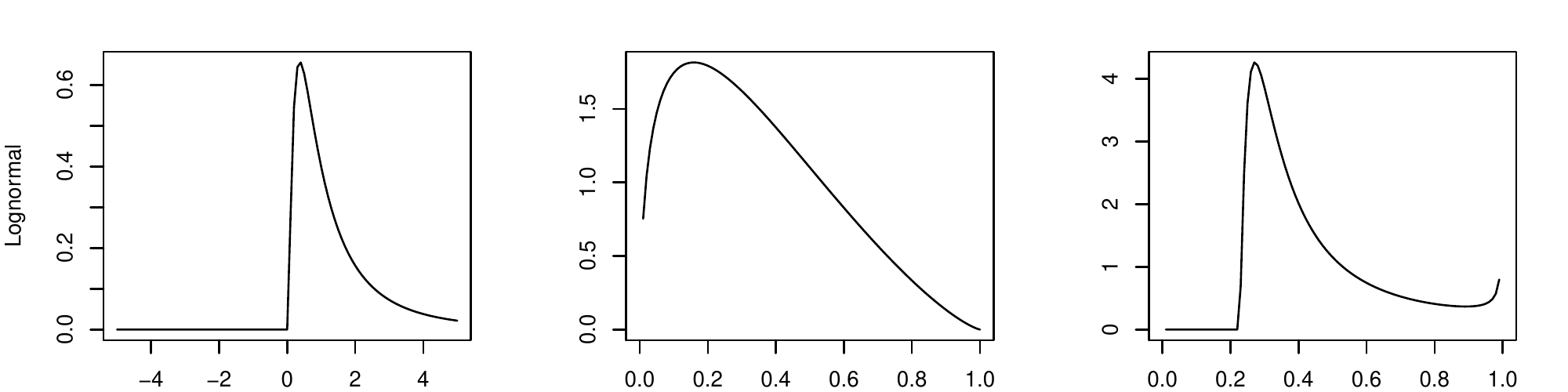}
\includegraphics[width=0.95\textwidth,height=0.95\textheight,keepaspectratio]{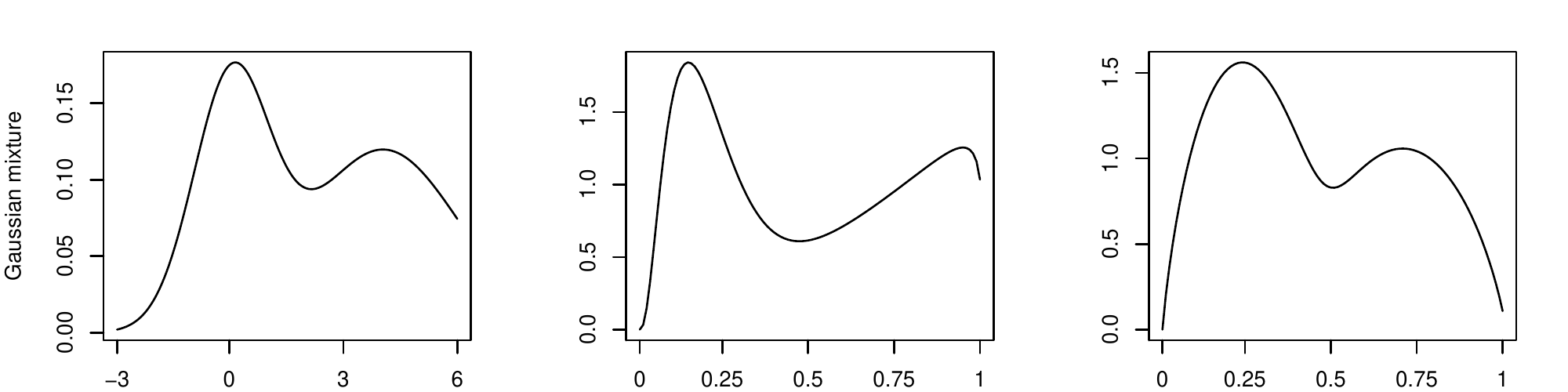}
\caption{Comparison of $f$, $f^*$ and $\tilde f$ for five distributions. }
\label{fig1}
\end{figure}

We next provide examples on the probability density functions $f$, $f^*$ and $\tilde f$
when $f$ is the density of Gaussian, Laplace, Lognormal and Uniform distributions. Also  a mixture of two Gaussian distributions denoted by $GMM(\mu_1,\mu_2,\sigma_1,\sigma_2,w)$
is considered with the densities
$
w\varphi_{\mu_1,\sigma_1}(x)+(1-w) \varphi_{\mu_2,\sigma_2}(x)
$,
$0\le w\le 1$.
Figure~\ref{fig1} then shows the impact of the transformations $f\to f^*$ and $f\to \tilde f$ in these cases.

\begin{longtable}[hth!]{lrrr>{\hspace*{.001\linewidth}}c<{\hspace*{.001\linewidth}}rrr}
\toprule
Distribution & \(e^{2H(f)}\) & \(e^{2H(f^*)}\) &
\(e^{2H(\tilde{f})}\) && \(H^*(f)^{-2}\) & \(H^*(f^*)^{-2}\) &
\(H^*(\tilde{f})^{-2}\)\tabularnewline
\midrule
\endhead
N(0,1) & 17.079 & 0.824 & 1.000 && 12.566 & 0.750 & 1.000\tabularnewline
Laplace(1) & 29.556 & 0.680 & 0.887 && 16.000 & 0.719 & 0.783\tabularnewline
Lognormal(0,1) & 17.079 & 0.642 & 0.308 && 7.622 & 0.537 & 0.186\tabularnewline
U(0,1) & 1.000 & 1.000 & 0.703 && 1.000 & 1.000 & 0.567\tabularnewline
GMM($0,4,1,2,0.4$) & 100.000 & 0.862 & 0.855 && 78.000 & 0.792 & 0.756\tabularnewline
\bottomrule
\caption{The power entropy and the \([H^*]^{-2}\) measure for some continuous distributions and their transformations.}
\label{t1}
\vspace{-3mm}
\end{longtable}

Table~\ref{t1} provides for the same distributions the power entropies $e^{2H(\cdot)}$ and $H^*(\cdot)^{-2}$ for  $f$, $f^*$ and $\tilde f$.
Note that the information measures applied to $f$ are not invariant under rescaling of $x$ as opposed to
$f^*$ and $\tilde f$. For example for the settings we use in the Table~\ref{t1}, the  normal and lognormal densities have the same power entropy just by accident and the equality is not generally true.

\begin{figure}[h!]
\begin{minipage}{1\textwidth}
\subcaption{GMM($0,\mu_2,1,1,0.5$) where $\mu_2$ varies.}
\vspace*{2pt}
\includegraphics[width=0.95\textwidth,height=0.95\textheight,keepaspectratio]{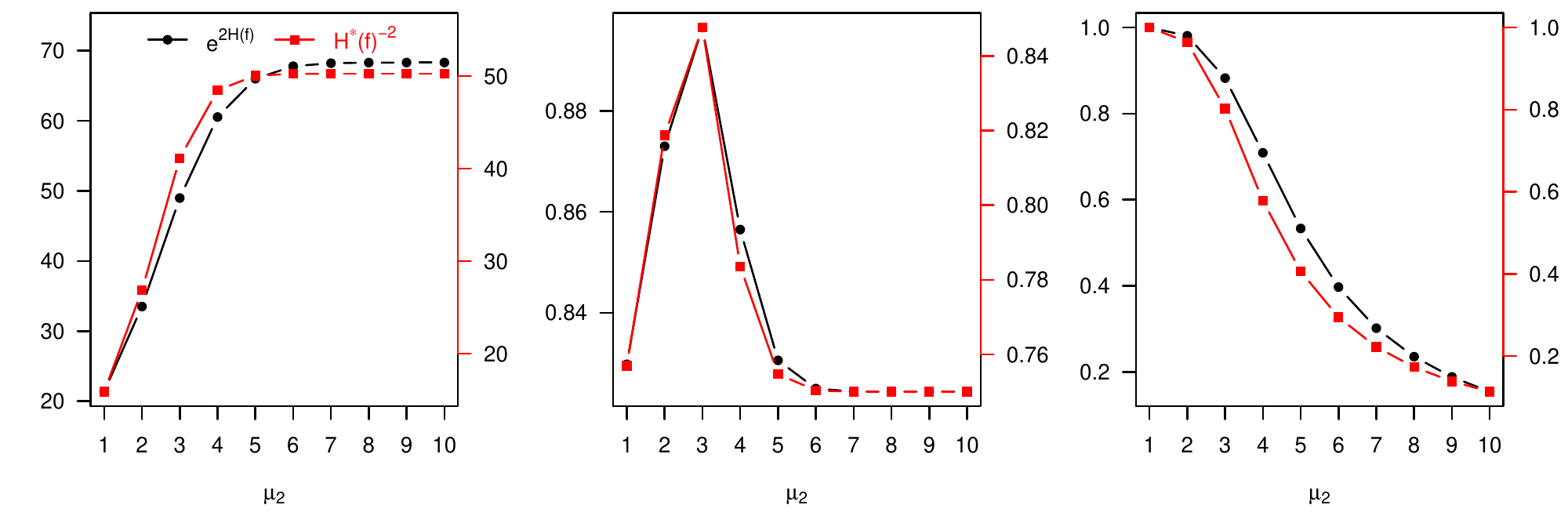}

\end{minipage}
\begin{minipage}{1\textwidth}
\subcaption{GMM($0,2,1,\sigma_2,0.5$) where $\sigma_2$ varies.}
\includegraphics[width=0.95\textwidth,height=0.95\textheight,keepaspectratio]{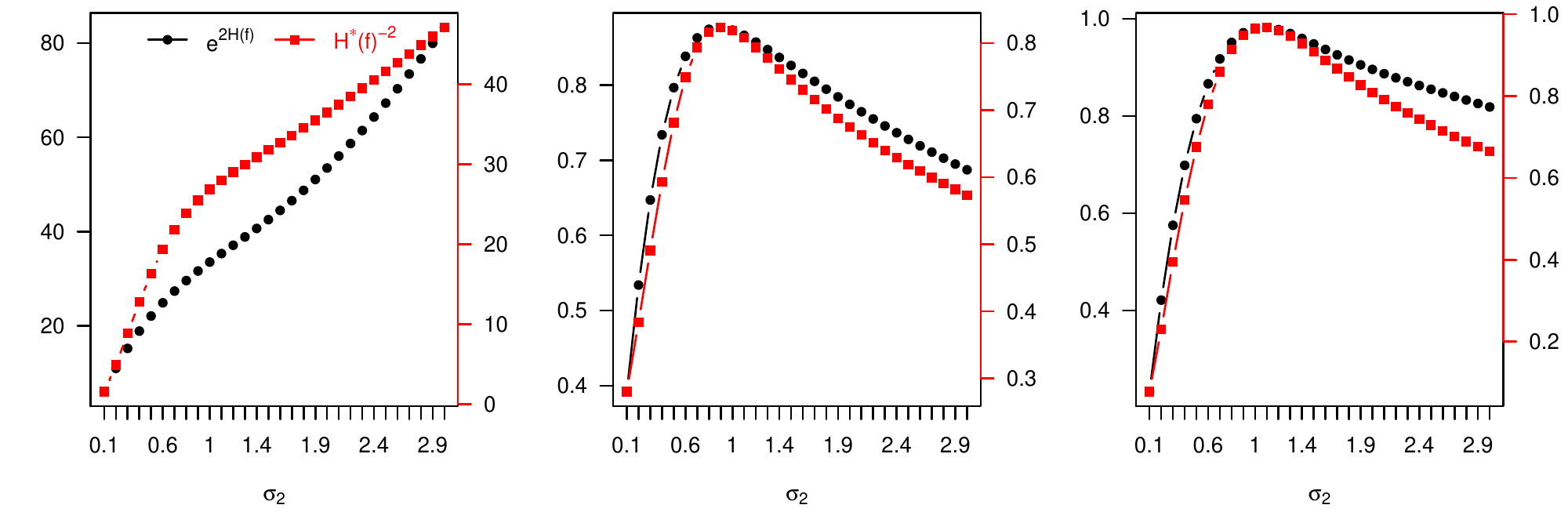}

\end{minipage}
\begin{minipage}{1\textwidth}
\subcaption{GMM($0,2,1,1,w$) where $w$ varies.} \label{fig:2c}
\includegraphics[width=0.95\textwidth,height=0.95\textheight,keepaspectratio]{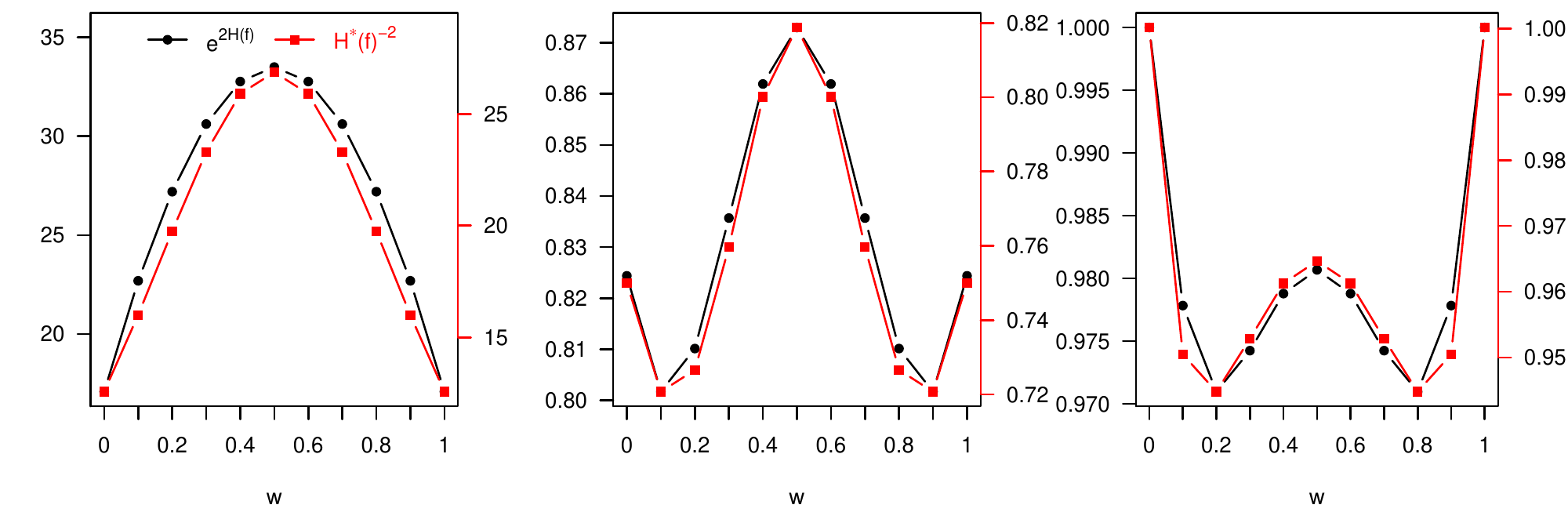}

\end{minipage}
\caption{Power entropy and \([H^*]^{-2}\) for different GMMs when always one parameter varies. The left vertical axis corresponds to power entropy and the right axis to \([H^*]^{-2}\). The left panel gives the measures for $f$, the middle for $f^*$ and the right for $\tilde f$.}
\label{fig2}
\end{figure}

For better understanding  on the measures, we illustrate the behavior of  $e^{H(\cdot)}$ and $H^*(\cdot)^{-2}$ in the GMM model with four fixed and one varying parameter, each in turn.
In Figure~\ref{fig2}  both information measure curves  are plotted in the same figures to compare the shapes of the curves as well as the occurrences of extreme values. The curves for $\tilde f$ with
varying location and scale seem natural as minimum information is attained as GMM gets ``closer'' to the normal distribution. Results for $f^*$ and varying location seem strange in a sense where one would expect decreasing behaviour of both measures as the distance in means increases, as it is case for $\tilde{f}$, while the result for $f$ in all three cases could simply be explained with decrease in information as a result of increase in overall variance of the mixture.
$e^{H(\cdot)}$ and $H^*(\cdot)^{-2}$  seem to behave almost proportionally in all cases. In cases of $f^*$ and $\tilde{f}$ where the majorization is well defined, such behaviour is indeed expected, as the reciprocals of both $e^{H(\cdot)}$ and $H^*(\cdot)^{-2}$ are information measures for both $f^*$ and $\tilde{f}$. However, further investigations into this matter will be conducted in the future.

\section{Independent component analysis}\label{sec3}

\subsection{Some preliminaries}\label{subsec3.1}
In this section we consider multivariate random variables.
For  a $p$-variate random vector $\bo x$ with finite second moments, the mean vector and covariance matrix are  $E(\bo x)\in \mathbb{R}^p$ and  $Cov(\bo x)\in \mathbb{R}^{p\times p}$, respectively.
Let $Cov(\bo x)=\bo U\bo D \bo U'$ be the eigenvector-eigenvalue decomposition of the covariance matrix. Then $Cov(\bo x)^{-1/2}:=\bo U\bo D^{-1/2} \bo U'$  and
$\bo x^{st}= Cov(\bo x)^{-1/2}(\bo x-E(\bo x))$
standardizes $\bo x$, that is, $E(\bo x^{st})=\bo 0$ and  $Cov(\bo x^{st})=\bo I_p$.
The set of $p\times r$, $r\le p$,  matrices with orthonormal columns is denoted by $\mathcal{O}^{p\times r}$. Thus
$
\bo U\in  \mathcal{O}^{p\times r}$ implies $\bo U'\bo U=\bo I_r$.
The set of $p\times p$ diagonal matrices with positive diagonal elements  is denoted by $\mathcal{D}^{p\times p}$.
If $\bo U\in  \mathcal{O}^{p\times p}$  and $\bo D\in \mathcal{D}^{p\times p}$  then
$\bo x\to \bo U \bo x$
and
$\bo x\to \bo D \bo x$, $\bo x\in \mathbb{R}^p$,
are a {\it rotation} operator  and a {\it componentwise rescaling}
operator, respectively. Let $\bo A\in \mathbb{R}^{p\times q}$ be a matrix with rank $r\le \min\{p,q\}$. Then the linear operator $\bo A$ may be written as (singular value decomposition, SVD)
$
\bo A=\bo U \bo D \bo V'=\sum_{i=1}^r d_i \bo u_i\bo v_i'
$
where $\bo U=(\bo u_1,...,\bo u_r)\in \mathcal{O}^{p\times r}$, $\bo V=(\bo v_1,...,\bo v_r)\in \mathcal{O}^{q\times r}$, and $\bo D\in \mathcal{D}^{r\times r}$.

\subsection{Elliptical model and independent components model}\label{subsec3.2}

 Let $\bo x$ be a $p$-variate vector with the full-rank covariance matrix $Cov(\bo x)$.  We say that $\bo x$ has a {\it spherical distribution} if there exists a $\bs \mu$ such that $(\bo x-\bs\mu)\sim \bo U(\bo x- \bs\mu)$ for all orthogonal $\bo U$. In the following we first define the elliptic and independent components distributions (see for example \cite{IC2,IC4} for more details).

\begin{Def}
Let $\bo x \in \mathbb{R}^p$ be a $p$-variate random vector.
\begin{enumerate}
\item $\bo x$ has an {\it elliptical distribution} if there exists
a nonsingular $\bo A\in \mathbb{R}^{p\times p}$ such that $\bo A\bo x$ has a spherical distribution. \
\item $\bo x$ has an {\it independent components distribution} if there exists
a nonsingular $\bo A\in \mathbb{R}^{p\times p}$ such that $\bo A\bo x$ has independent components.
\end{enumerate}
\end{Def}

We next provide some results on how the matrix $\bo A$ can be found in different cases.

\begin{The}\label{T3.1}
Let $\bo x$ be a $p$-variate random vector with a full-rank covariance matrix $Cov(\bo x)=\bo U\bo D \bo U'$. 
Then we have the following.
\begin{enumerate}
\item  $[\bo V \bo D^{-1/2}\ \bo U']\bo x$ has uncorrelated components for all orthogonal $\bo V$.
\item If $\bo x$ has an elliptical distribution,  $[\bo V \bo D^{-1/2}\ \bo U']\bo x$ has a spherical distribution for all orthogonal $\bo V$.
\
\item If $\bo x$ has an independent components distribution,  $[\bo V \bo D^{-1/2}\ \bo U']\bo x$ has independent components for some choice(s) of orthogonal $\bo V$.
\item  If $\bo x$ has both an elliptical distribution and an independent component distribution
then  $[\bo V \bo D^{-1/2}\ \bo U']\bo x$ has independent Gaussian components for all  orthogonal $\bo V$, that is, $\bo x$ has a multivariate Gaussian distribution.
\end{enumerate}
\end{The}

\subsection{Projection pursuit and independent component analysis}\label{subsec3.3}

 Let  
$\bo x$ have an independent components distribution such that $\bo z=\bo A\bo x+\bo b$ is standardized ($E(\bo z)=\bo 0$ and $Cov(\bo z)=\bo I_p$) and has independent components.
Theorem~\ref{T3.1} then implies that $\bo A=\bo V' Cov(\bo x)^{-1/2}$ where the rotation matrix
$\bo V$ can be chosen as $\bo V=(\bo V_1,\bo V_2)$ separating non-Gaussian independent components in
$\bo V_1' Cov(\bo x)^{-1/2}\bo x$ and Gaussian independent components in $\bo V_2' Cov(\bo x)^{-1/2}\bo x$. Note that $\bo V_2$ is only unique up to right multiplication by an orthogonal matrix.
A generally accepted strategy  is to find
$\bo V_1=(\bo v_1,...,\bo v_q)\in \mathcal{O}^{p\times q}$
such that the components of $\bo V_1'\bo x^{st}$ are as `non-Gaussian as possible'.
The Gaussian part $\bo V_2' Cov(\bo x)^{-1/2}\bo x$ is thought to be just the noise part and, for other components,
it is argued that the sum of independent random variables is `more Gaussian' than the original variables.
The noise interpretation of the Gaussian part may be motivated by the following.
A random vector  has a multivariate normal distribution if and only if all linear combinations of the marginal variables have univariate normal distributions, that is, there are no `interesting' directions. The normal distribution is the only distribution for which all third and higher cumulants are zero. As seen before, a Gaussian distribution is the  distribution with the poorest information among distributions
with the same variance (highest entropy, smallest Fisher information). \
For a thorough discussion of Gaussian distributions, see \cite{Shevlyakov}.

Let $D(x)$ then be the {\it projection index}, i.e., the functional that is used to measure non-Gaussianity. In the one-by-one projection pursuit approach the first direction $\bo v_1$  ($\bo v_1'\bo v_1=1$) maximizes $D(\bo v_1'\bo x^{st})$, the second
direction $\bo v_2$ is orthogonal to $\bo v_1$ ($\bo v_2'\bo v_2=1, \bo v_2'\bo v_1=0$) and maximizes $D(\bo v_2'\bo x^{st})$ and so on. After finding $\bo v_1,...,\bo v_{j-1}$, we optimize
the Lagrangian function
\[ L(\bo v;\lambda_{j1},...,\lambda_{jj})=  D(\bo v'\bo x^{st}) 
-\lambda_{jj} (\bo v' \bo v-1)-\sum_{i=1}^{j-1} \lambda_{ji} \bo v' \bo v_i.\]
Then $\bo v_j$ solves the (estimating) equation
$
( \bo I_p-\sum_{i=1}^{j-1} \bo v_i \bo v_i') \bo T(\bo v)=(\bo T(\bo v)'\bo v) \bo v,
$
where
$
\bo T(\bo v)=  \partial D(\bo v'\bo x^{st})/\partial \bo v  
$.
From the computational point of view, this suggests a {\it fixed-point  algorithm}. The estimation equation
also provides a way to find the limiting distribution of the estimate, since the estimate is obtained when
the theoretical multivariate distribution is replaced by the empirical one.
See for example \cite{IC5,IC1,IC6} and references therein for more details.

The following questions naturally arise.  How should one choose the projection index $D(x)$ to find the  independent components? Are the independent components provided by the  most informative directions as has been often stated in the literature? These questions are partially answered by the following.

\begin{The}\label{T3.2}
Let $\bo z=\bo A\bo x+\bo b=(z_1,...,z_p)'$ be the vector of standardized
independent components.
\begin{enumerate}
\item Let $D(x)$ be a subadditive squared dispersion measure. \\ Then $D(\bo v'\bo x^{st} )\le \max_j D(z_j)$.
\item Let $D(x)$ be a superadditive squared dispersion measure.\\ Then $D(\bo v'\bo x^{st})\ge  \min_j D(z_j)$.
\end{enumerate}
\end{The}

Based on Theorem~\ref{T3.2} and the discussion above
we can now end the paper with the following conclusions.
If $D(x)$ is subadditive then it can be used as a projection index. For example the cumulants
$\kappa_{2k+1}^{2/(2k+1)}(x)$ and  $\kappa_{2k+2}^{2/(2k+2)}(x)$, $k=1,2,\ldots$,
when calculated for standardized distributions,  provide squared dispersion measures that are subadditive.  Therefore they can be used as projection indices.
For superadditive $D(x)$, the functional $(D(x))^{-1}$ is a valid projection index
as  $(D(\bo v'\bo x^{st} ))^{-1}\le \max_j (D(z_j))^{-1} $.
As seen before, entropy power $e^{H(x)}$ and the inverse of Fisher information, $J^{-1}(x)$ are superadditive squared dispersion measures.
Note that in both cases  $D(\bo v'\bo x^{st})$ is in fact a ratio of two squared dispersion functions, and the projection index measures deviation from Gaussianity using a skewness,  kurtosis or information measure. As mentioned in Section~\ref{subsec3.3}, $(\kappa_3^2(x)+(1/4) \kappa_4^2(x))/12$ provides
an approximation of negative differential entropy in a neighborhood of Gaussian distribution and is a valid projection index as well. For further discussion, see \cite{huber}.
Note also that one of the most popular ICA procedures in the engineering community,  the so called {\it fastICA},  uses a projection index of the form  $D(x)=|E[C(x)])|$ where $C$ is such a function that, if $z\sim N(0,1)$ then  $E[C(z)]=0$. Examples of valid choices of $C$ are
$C(z)=z^3$ and $C(z)=z^4-3$ providing again the third and fourth cumulants, respectively.

\section{Final remarks}\label{sec4}

The usage of various information criteria is popular in  independent component analysis.
The connections between notions of information and statistical independence and the special role of the Gaussian distribution were discussed in detail in the paper.
We also introduced new ideas and  partial orderings for information which utilize  transformed location and scale-free probability density functions. In  independent component analysis with unknown marginal densities, the estimation of the value of the adapted information measure in a given direction is highly challenging
and it has to be done again and again  when applying the fixed point algorithm for the correct direction.
Substantial research is therefore still needed for these tools to be of practical value.

\section{Appendix: The Proofs}\label{Appendix}

{\bf Proof of Theorem~\ref{additivity theorem}.}
Let $x_1,...,x_n$ be a random sample from a distribution of $x$ with the mean value $E(x)$ and variance $Var(x)$. By the central limit theorem,
\[
z_n=\frac 1{\sqrt{n}} \sum_{i=1}^n \frac{x_i-E(x)} {\sqrt{Var(x)}}\to_d z\sim N(0,1).
\]
 Therefore, by additivity and affine equivariance,
\[
T(z_n)=\sqrt{\frac n{Var(x)}}(T(x)-E(x))\to 0\ \ \mbox{and}\ \
S^2(z_n)=\frac{S^2(x)}{Var(x)}\to S^2(z)
\]
and the result follows. For similar results in the multivariate case, see \cite{virta2018characterizations}.
\\

{\bf Proof of Theorem~\ref{dispersion measures}.}
The invariances of the measures $H(x)$, $H^*(x)$, $H^{**}(x)$ and $J(x)$ under location
shifts $f(x)\to f(x+b)$, sign change $f(x)\to f(-x)$ as well as under $f\to  f_{a,b,\Delta}$ follow easily from their definitions and from the definition of the Riemann integral.
We therefore only have to consider rescaling $f(x)\to (1/a)f(x/a)$ with $a>0$. Then $H(ax)=-\int (1/a)f(x/a)\log((1/a)f(x/a))dx=-\int f(x)\log((1/a)f(x))dx
=H(x)+\log(a)$ and therefore $e^{2H(ax)}=a^2e^{H(x)}$. In a similar way one can show that $[H^*(ax)]^{-2}=a^2 [H^*(x)]^{-2}$. Also easily $[H^{**}(ax)]^{-2}=a^2 [H^{**}(x)]^{-2}$. As $f'(x)\to (1/a^2)f'(x/a)$
one also easily shows that $[J(ax)]^{-1}=a^2 [J(x)]^{-1}$. Thus all the four measure are scale equivariant and therefore squared dispersion measures.
\\

{\bf Proof of Theorem~\ref{KL}.} $f:g$ is indeed a density function since it is trivially nonnegative and $\int_0^1 (f:g)(u)\mathrm{d}u=\int_0^1 \frac{f(G^{-1}(u))}{g(G^{-1}(u))}\mathrm{d}u=\int_{-\infty}^\infty f(x)\mathrm{d}x=1$
with the substitution $x=G^{-1}(u)$. Similary, $-H(f:g)=\int_0^1 (f:g)(u)\log((f:g)(u))\mathrm{d}u
=\int_{-\infty}^{\infty} f(x)\log\frac{f(x)}{g(x)}\mathrm{d}x=D(f||g).
$
\\

{\bf Proof of Theorem~\ref{T3.1}.}
(1) Let $\bo V$ be orthogonal.  As $Cov([\bo V \bo D^{-1/2} \bo U']\bo x)=\bo V \bo D^{-1/2}\bo U' Cov(\bo x)\bo U \bo D^{-1/2} \bo V'=\bo V\bo V'=\bo I_p$, the components  of $[\bo V \bo D^{-1/2} \bo U']\bo x)$ are uncorrelated.
(2) Assume that $\bo A\bo x$ is spherical with $\bo A=\bo V\bo C\bo W'$ rescaled so that $Cov(\bo A\bo x)=\bo I_p$. As $\bo A Cov(\bo x)\bo A'=\bo I_p$, $Cov(\bo x)=(\bo A'\bo A)^{-1}$ and  $\bo W \bo C^{-2} \bo W'=\bo U \bo D \bo U'$. Therefore $\bo W=\bo U$ and $\bo C=\bo D^{-1/2}$ and we can conclude that $[\bo V\bo D^{-1/2}\bo U']\bo x$ is spherical for any orthogonal $\bo V$. (If $\bo x$ is spherical then
$\bo V\bo x$ is spherical for all orthogonal $\bo V$.)
(3) Let $\bo A\bo x$  with $\bo A=\bo V\bo C\bo W'$ have independent and standardized components so that $Cov(\bo A\bo x)=\bo I_p$. As in (2), $\bo A$ must be $\bo V\bo D^{-1/2}\bo U'$ but now for some $\bo V$ only. (It is not true that if $\bo x$ has independent standardized components then $\bo V\bo x$ has independent components for any choice of  $\bo V$.)
(4) Based on (2) and (3), there exist an $\bo A=\bo V\bo D^{-1/2}\bo U'$ such that $\bo A\bo x$
has a spherical distribution with independent components. Then by the Maxwell-Hershell theorem,
$\bo A\bo x$ has a multivariate normal distribution. For the proof of the Maxwell-Hershell theorem,
see e.g. Proposition 4.11. in  \cite{BilodeauBrenner1999}.

{\bf Proof of Theorem~\ref{T3.2}.} Let $\bo z=\bo A\bo x+\bo b=(z_1,...,z_p)'$ be a vector of standardized independent components. By Theorem~\ref{T3.1}, $\bo z=\bo V\bo x^{st}$ with some orthogonal $\bo V$.
If $\bo u'\bo u=1$ then also $(\bo V\bo u)'(\bo V\bo u)=1$ and therefore
$D(\bo u'\bo x^{st})=D(\bo u'\bo V\bo z)\leq \sum (\bo V'\bo u)_i^2 D(z_i)\leq \max_j D(z_j)$
for subadditive squared dispersion measure $D$ and
$D(\bo u'\bo x^{st})=D(\bo u'\bo V\bo z)\geq \sum (\bo V'\bo u)_i^2 D(z_i)\geq \min_j D(z_j)$
 for superadditive squared dispersion measure $D$.

\section{Acknowledges}
  The work of KN has been supported by the Austrian Science Fund (FWF) Grant number P31881-N32.

\end{document}